\documentclass{amsart}
\usepackage{amssymb,latexsym}

\newtheorem{thm}{Theorem}[section]
\newtheorem{cor}[thm]{Corollary}
\newtheorem{lem}[thm]{Lemma}

\newtheorem{rem}[thm]{Remark}

\newcommand{\Rn}{\mathbb{R}^{n}}

\title{Lemma Poincar\'e for $L_{\infty,\,loc}$-forms}

\author{Vladimir Gol'dshtein, Stanislav Dubrovskiy}

\address{Department of Mathematics, Ben Gurion
University of the Negev, P.O.B. 653\\Beer-Sheva 84105, Israel}

\email{vladimir@math.bgu.ac.il, dubr@math.bgu.ac.il}

\keywords{Lipschitz manifold, $L_p$-cohomology, Lemma Poincar\'e}

\date{December 11, 2007}

\begin{document}

\maketitle

\begin{abstract}
We show that every closed $L_{\infty,\,loc}$-form on
$\mathbb{R}^{n}$ is exact.\\Differential is understood in the
sense of currents. The proof does not use any explicit geometric
constructions. De Rham theorem follows.
\end{abstract}


\section{Introduction}
We work with a class $\Omega^{*}_{\infty,\,loc}$ of differential
forms on $\mathbb{R}^{n}$, those with coefficients locally in
$L_{\infty}$, with differential satisfying the same condition.
``Locally'' will mean ``on a neighborhood with compact closure'',
as usual. Classes $\Omega^{*}_{\infty}$ and
$\Omega^{*}_{\infty,\,loc}$ are invariant under bi-Lipschitz
coordinate changes and therefore are naturally defined on any
Lipschitz manifold. In context of his geometric theory of
integration Whitney in his book \cite{W} considered a class of
so-called flat forms, invariant under Lipschitz maps, and proved a
Poincar\'e lemma for them \cite[Th IX.12A]{W} using an elaborate
geometric construction.  
We give a different shorter proof below, based on the techniques
from \cite{GKS}. 
Differential forms of class $\Omega^{*}_{\infty,\,loc}$ are
nothing but locally flat forms in the sense of Whitney.

The second author is grateful to the first author for introducing
him to the subject and would also like to acknowledge the support
of the Skirball foundation and hospitality of the Center of
Advanced Studies in Mathematics at Ben Gurion University.

\section{Class $\Omega^{*}_{\infty,\,loc}$}
The differential forms we consider are not the ones in the smooth
category, they belong properly in the Lipschitz category and are
defined as follows.

On an open domain $U\subset\mathbb{R}^{n}$ a $k$-form is an almost
everywhere defined map
$$w:U\rightarrow\Lambda^{k}(\mathbb{R}^{n})$$
Two maps are considered equivalent and define the same form if
they coincide almost everywhere. Denote by $\mathcal{F}^{k}(U)$ a
vector space of differential forms of degree $k$, defined on $U$.
Vector spaces $\mathcal{F}^{k}(U)$ form a graded ring,
$\mathcal{F}^{*}(U)$ with respect to addition and exterior
product.

By $\mathcal{L}^{k}_{\infty}$ we denote the space of $k$-forms
with bounded coefficients. The union of these for all $k$ make up
a graded ring $\mathcal{L}^{*}_{\infty}$.

The ring of smooth forms with compact support in $U$ will be
denoted $\mathcal{D}^{*}(U)$. According to de Rham \cite{dR},
\emph{$k$-currents} on $U$  are linear continuous functionals on
$\mathcal{D}^{n-k}(U)$, denoted $\mathcal{E}^{k}(U)$. For each
current $w\in\mathcal{E}^{k}(U)$ one can define its exterior
differential $dw\in\mathcal{E}^{k+1}(U)$ by the formula:
$$<dw,\alpha>=(-1)^{k+1}<w,d\alpha>,\quad\quad \forall\alpha\in\mathcal{D}^{n-k-1}(U)
$$
Define
$$\Omega^{*}_{\infty}(U)=\{w\mid w\in\mathcal{L}^{*}_{\infty}(U),\,dw\in\mathcal{L}^{*}_{\infty}(U)\}\ ,$$
where by $dw\in\mathcal{L}^{*}_{\infty}(U)$ we mean that $dw$ as a
functional has an integral representation with an
$\mathcal{L}^{*}_{\infty}$-kernel:
$$<dw,\alpha>=\int_{U}dw\wedge\alpha,
$$
where we abuse notation: the $dw$ under the integral means the
form $dw\in\mathcal{L}^{*}_{\infty}(U)$, same notation as the
original functional it represents. Using a multi-index
$$\vec{i}=(i_{1},i_{2},...,i_{k})\quad\quad i_1<...<i_{k}\,,$$ we shall write
a $k$-form $w$ as follows:
$w=\displaystyle{\sum_{\vec{i}}a_{\vec{i}}dx^{\vec{i}}}$.\\On
every compact $K$ in $U$ consider a seminorm:
$$|w|_{K,\infty}=\max_{\vec{i}} \textrm{ess}\sup_{x \in K}|a_{\vec{i}}(x)|
$$
We will use $\sup$ to denote essential supremum from now on. In
$\Omega^{*}_{\infty}$ consider\\seminorms
$$|w|_{K,\infty,\infty}=|w|_{K,\infty}+|dw|_{K,\infty}
$$
With these seminorms $\mathcal{L}^{k}_{\infty}(U)$ and
$\Omega^{k}_{\infty}(U)$ become locally convex
topological vector spaces. 

We will compare now the forms of class $\Omega^{*}_{\infty,\,loc}$
with locally flat forms in the sense of Whitney \cite[IX.6]{W}. We
recall some definitions.

If $P$ is an $s$-dimensional plane in $\Rn$, and $A$ - a
measurable set in $P$, then $|A|_s$ shall stand for the
$s$-dimensional Lebesque measure of $A$.

Given a form $w\in\mathcal{F}^{k}(U)$ and a set $Q\subset U$, a
simplex $\sigma^s\subset U$ $(s\geq k)\ $ is said to be
$Q$-\emph{good} for $w$, if all coefficients of $w$ are measurable
on $\sigma^s$ and $|\sigma^s\backslash Q|_s=0$. Simplex $\sigma^s$
is said to be $Q$-\emph{excellent} if $\sigma^s$ and each of its
faces of dimension $\geq k$ is $Q$-good.

A form $w=\displaystyle{\sum_{\vec{i}}a_{\vec{i}}dx^{\vec{i}}}\in
\mathcal{F}^{k}(U)$ is said to be \emph{flat} if:

(a) coefficients $a_{\vec{i}}$ are measurable in $U$,

(b) $\displaystyle{\sup_{x\in U}|a_{\vec{i}}(x)|\leq N}$ for some
$N$,

(c) there exist such $Q\subset U$ and $N'$, that $|U\backslash
Q|_n=0$, and for every\\ $Q$-excellent for $w$ simplex
$\sigma^{k+1}\subset U$: $$|\int_{\partial \sigma^{k+1}}w|\leq
N'|\sigma^{k+1}|_{k+1}\ .$$
\begin{thm}
\label{thm:1.5} Class $\Omega^{*}_{\infty,\,loc}(U)$ coincides
with the space of all locally flat forms in $U$. The exterior
differential of an arbitrary form $w\in
\Omega^{*}_{\infty,\,loc}(U)$ coincides with the differential of
$w$ in the sense of Whitney \cite[IX.12]{W}.
\end{thm}
{\bf Proof}\quad See \cite[Theorem 1.5]{GKS} $\Box$

\section{$\Omega^{*}_{\infty,\,loc}$-morphisms}
A map $f:U\rightarrow \mathbb{R}^{m}$, $U\subset \Rn$, is called
locally Lipschitz if for every compact $K \subset U$ there exists
a constant $C$, satisfying: $$|f(x)-f(x')|\leq C|x-x'|$$ for any
$x, x' \in K$ . Every locally Lipschitz function $f:U\rightarrow
\mathbb{R}^{1}$ is differentiable almost everywhere on $U$, its
partial derivatives are generalized partial derivatives. Hence,
the exterior differential $df$ of a locally Lipschitz function
$f$, considered as a form of degree $0$, has a coordinate
representation $$df=\frac{\partial f}{\partial
x^1}dx^1+...+\frac{\partial f}{\partial x^n}dx^n\ .$$

Locally Lipschitz map $\varphi: U\rightarrow V$ will be called an
$l$-map, if for any point $x\in U$ there exist a neighborhood $W$
of $x$ and a constant $C$, such that $$|\varphi^{-1}(A)\cap
W|_n\leq C|A|_m$$ for any open $A\subset V$. Preimage of a measure
zero set under an $l$-map has measure zero. Hence for any $l$-map
$\varphi: U\rightarrow V$ there is a well-defined homomorphism
$\varphi^*: \mathcal{F}^*(V)\rightarrow \mathcal{F}^*(U)$.

According to Rademacher theorem \cite[IX.11]{W} any locally
Lipschitz map is differentiable almost everywhere, therefore any
$l$-map $\varphi$ satisfies:
$$\varphi^*(\mathcal{L}^{*}_\infty(V))\subset
(\mathcal{L}^{*}_\infty(U))\ .$$
\begin{thm}
\label{thm:2.2} If $\,\varphi: U\rightarrow V$, $U\subset \Rn$,
$V\subset \mathbb{R}^{m}$ - an $\,l$-map, $w\in
\Omega^{*}_{\infty,\,loc}(V)$, then $$\varphi^* w\in
\Omega^{*}_{\infty,\,loc}(U) ,\quad d(\varphi^* w)= \varphi^* dw\
.$$ The map
$$\varphi^*: \Omega^{*}_{\infty,\,loc}(V) \rightarrow \Omega^{*}_{\infty,\,loc}(U)$$
is continuous.
\end{thm}
{\bf Proof}\quad See \cite[Theorem 2.2]{GKS} $\Box$
\begin{rem}
A homeomorphism $\varphi: U \rightarrow V\,,\,U,\,V \subset \Rn $
is called Lipschitz if $\varphi$, $\varphi^{-1}$ are locally
Lipschitz. For Lipschitz homeomorphisms Theorem \ref{thm:2.2} is
satisfied, since they are $l$-maps.\end{rem} Whitney defines
homomorphism $\varphi^*$ on classes of flat forms for any locally
Lipschitz $\varphi: U \rightarrow V\,,\,U\subset \Rn,\,V  \subset
\mathbb{R}^{m}$ \cite[X.9]{W} . Due to Theorem \ref{thm:1.5} the
construction of Whitney defines a homomorphism $\varphi^*:
\Omega^{*}_{\infty,\,loc}(V) \rightarrow
\Omega^{*}_{\infty,\,loc}(U)\,$. If $\varphi$ is a Lipschitz
homeomorphism then $\varphi^*$ in the sense of Whitney
\cite[X.9]{W} coincides with the $\varphi^*$ introduced above.

\section{Lemma Poincar\'e for bounded forms}
Lemma Poincar\'e for forms with $p$-summable coefficients for
finite $p$ is proven in \cite[Lemma 3.4]{GKS}. The case $p=\infty$
(locally flat forms) is treated in \cite{W} via nontrivial
explicit geometric construction. We use analytic technique similar
to that used in \cite{GKS} to obtain a short alternative proof for
$p=\infty$ case. We start with the following key observation.
\begin{lem}
\label{lem:3.2} Let $I$ be a cube in $\mathbb{R}^{n}$ and $w$ an
arbitrary smooth $k$-form on $I$. There exists a smooth closed
$w'\ (dw'=0)$, such that:
$$|w-w'|_{I,\infty}\leq
C|dw|_{I,\infty}$$ for some constant $C$ depending only on $n$,
$k$ and $I$.
\end{lem}
{\bf Proof}\quad We can assume $I=I^{n}$ in $\mathbb{R}^{n}$. For
$k\geq n$ the statement is true for $C=0$, for other cases we will
proceed by induction on $n$. For $n=1$ we only need to consider
$k=0$. In this case $w=f(x)$ and we have:
$$\sup_{x\in I}|f(x)-f(0)|=\sup_{x\in I}|\int_{0}^{x}f'(\tau)d\tau|\leq\int_{0}^{1}|f'(\tau)|d\tau\leq\sup_{x\in
I}|f'(x)|\ .$$ In other words, $$|f(x)-f(0)|_{I,\infty}\leq
|df|_{I,\infty}\ .$$ Suppose now $n>1$. An arbitrary smooth form
$w$ can be represented in the following way:
$$w=w_{1}\wedge dx^{n}+w_{2}\ ,$$
where forms $w_1$ and $w_2$ do not contain $dx^n$ in their
coordinate representations. Clearly:
\begin{equation}
\label{eq:9}
|w|_{I,\infty}=\max(|w_1|_{I,\infty},|w_2|_{I,\infty})
\end{equation}
Forms $w_1$ and $w_2$ can be considered as those on the cube
$I'=I^{n-1}$, depending on a parameter $x^n$. Differentials of
these in $I'$ we will denote by $d'w_1$ and $d'w_2$. Then
$$
dw=d'w_1\wedge dx^n + (-1)^k\frac{\partial w_2}{\partial
x^n}\wedge dx^n + d'w_2\ ,$$ where by $\frac{\partial
w_2}{\partial x^n}$ we mean a form, whose coefficients are the
partial derivatives of $w_2$ with respect to $x^n$. We have:
\begin{equation}
\label{eq:10} (dw)_1=d'w_1+(-1)^k\frac{\partial w_2}{\partial
x^n}\ ,\ (dw)_2=d'w_2
\end{equation}
Denoting $$d'w_2=\sum_{i_1<...<i_{k+1}<n}a_{\vec{i}}dx^{\vec{i}}$$
we have:
$$|d'w_2|_{I,\infty}=\max_{\vec{i}=(i_1,...,i_{k+1})}
\sup_{x \in I}|a_{\vec{i}}(x)|=\max_{\vec{i}}\sup_{x^n \in
[0,1]}\sup_{x \in I'}|a_{\vec{i}}(x)|=\sup_{x^n \in
[0,1]}|d'w_2|_{I',\infty}
$$
Hence for some $n^{\textrm{th}}$ coordinate value $\tau \in
[0,1]$:
\begin{equation}
\label{eq:11} |d'w_2(\tau)|_{I',\infty}\leq |d'w_2|_{I,\infty}
\end{equation}
Consider a form $$\bar{w}=\int_{\tau}^{x^n}w_1(\sigma)d\sigma\ ,$$
with the integration carried out over the $n^{\textrm{th}}$
coordinate. Its differential is
\begin{equation}
\label{eq:12}
d\bar{w}=(-1)^{k+1}w+\int_{\tau}^{x^n}(dw)_1dx^n+(-1)^k w_2(\tau)\
.
\end{equation}
The induction step gives us a closed form $\tilde{w}$ on $I'$ such
that:
\begin{equation}
\label{eq:13} |w_2(\tau)-\tilde{w}|_{I',\infty}\leq
C|d'w_2|_{I',\infty}\ .
\end{equation}
Let us define $w'=(-1)^{k-1}d\bar{w}+\tilde{w}$. Here $\tilde{w}$
is considered as a form on $I$, independent of $x^n$. According to
(\ref{eq:12})
\begin{equation}
\label{eq:14}
w-w'=(-1)^k\int_{\tau}^{x^n}(dw)_1dx^n+w_2(\tau)-\tilde{w}\ .
\end{equation}
Let us estimate the norm of
$\displaystyle{\int_{\tau}^{x^n}(dw)_1dx^n}$. Since the integral
is taken separately of each coefficient $a_{\vec{i}}$ of the form
$(dw)_1$ we have:
$$|\int_{\tau}^{x^n}(dw)_1dx^n|_{I,\infty}=\max_{\vec{i}}\sup_{I}|\int_{\tau}^{x^n}a_{\vec{i}}dx^n|$$
$$\leq \max_{\vec{i}}\sup_{I}\int_{0}^{1}|a_{\vec{i}}|dx^n\leq
\max_{\vec{i}}\sup_{I}\sup_{x^n}|a_{\vec{i}}|=|(dw)_1|_{I,\infty}
$$
Hence, from (\ref{eq:14}), (\ref{eq:13}), (\ref{eq:11}),
(\ref{eq:10}), (\ref{eq:9}):
$$
|w-w'|_{I,\infty}\leq |(dw)_1|_{I,\infty} +
C|(dw)_2|_{I,\infty}$$
$$\leq C'(|(dw)_1|_{I,\infty}
+|(dw)_2|_{I,\infty})\leq 2C'\max(|(dw)_1|_{I,\infty}
,|(dw)_2|_{I,\infty})\ ,
$$ where $C'=\max(1,C)$. \\Now notice that
$$|dw|_{I,\infty}=|(dw)_1\wedge dx^n +
(dw)_2|_{I,\infty}=\max(|(dw)_1|_{I,\infty},|(dw)_2|_{I,\infty})\
.
$$
Thus $|w-w'|_{I,\infty}\leq 2C'|dw|_{I,\infty}$ which concludes
the proof.
 $\Box$
\begin{cor}
\label{cor:3.3} For any smooth closed $k$-form $w$ ($k>0$) on any
cube $I \subset \mathbb{R}^{n}$ \\there exists such smooth
"primitive" $\theta$ ($d\theta=w$) that
$$|\theta|_{I,\infty}\leq C|w|_{I,\infty}
$$
for some constant $C$.\end{cor} {\bf Proof}\quad According to the
(usual) Lemma Poincar\'e for smooth forms, there exists a smooth
$\theta_1$, such that $d\theta_1=w$. Lemma \ref{lem:3.2} then
provides for existence of a smooth closed $\theta_2$ satisfying
$|\theta_1-\theta_2|_{I,\infty}\leq C|w|_{I,\infty}$. The form
$\theta=\theta_1-\theta_2$ is the desired primitive. $\Box$
\begin{lem}
\label{lem:3.4} Any $x \in \mathbb{R}^{n}$ has a neighborhood $V$
such that any closed\\ $\mathcal{L}^{k}_{\infty}$-form $w$ on $V$
is exact: for $k>0$ $\exists \theta \in
\mathcal{L}^{k-1}_{\infty}(V): d\theta=w$. \\If $k=0$, then $w$ is
a locally constant function.
\end{lem}
{\bf Proof}\quad A function which has all of its generalized
derivatives vanish is locally constant. For $k>0$ consider an open
cube $I$, centered at $x$. According to \cite[Lemma 1.3]{GKS} and
the Remark following it, it is possible to choose a sequence of
smooth closed $k$-forms that weakly converges to $w$.

Using Mazur's theorem \cite[Theorem V.1.2]{Y} we can construct
another sequence \{$w_s$\} (of finite convex combinations of the
elements of the original sequence) that will satisfy:
$$|w-w_s|_{I,\infty}\leq \frac{1}{2^s}\ ,$$ and $dw_s=0$ on $I$.

From Corollary \ref{cor:3.3} it follows that there exist such
forms $\theta_s$ on the cube $I$, that
$$
d\theta_1=w_1\ ,\ d\theta_s=w_s-w_{s-1}\ \textrm{ for }s>1\
\textrm{ and }\ |\theta_s|_{I,\infty}\leq
C|w_s-w_{s-1}|_{I,\infty}\ .
$$
Since $$|\theta_s|_{I,\infty} \leq C|w_s-w_{s-1}|_{I,\infty} \leq
C(|w-w_{s}|_{I,\infty}+|w-w_{s-1}|_{I,\infty})\leq
\frac{C}{2^s}+\frac{C}{2^{s-1}}\ ,$$ the series
$\displaystyle{\sum_{s=1}^{\infty}\theta_s}$ converges in
$\mathcal{L}^{k-1}_{\infty}(I)$. Let
$\theta=\displaystyle{\sum_{s=1}^{\infty}\theta_s}$. \\Since
$$|w-\sum_{s=1}^{N}d\theta_s|_{I,\infty}=|w-w_{N}|_{I,\infty}\leq \frac{1}{2^N}\ ,
$$
we see that $d\theta=w$, as required. $\Box$\\

Denote by $H^*_{\infty,\,loc}(M)$ the cohomology of the cochain
complex $\{\Omega^{*}_{\infty,\,loc},d\}$ on a manifold $M$.
\begin{thm}[De Rham Theorem]
\label{thm:dR} For any Lipschitz manifold $M$\\ the cohomology
groups $H^*_{\infty,\,loc}(M)$ are canonically isomorphic to the
cohomology groups $H^*(M;\mathbb{R})$ of the manifold $M$.
\end{thm}
Denote by $\lambda : H^{*}_{\infty,\,loc}(M)\rightarrow
H^*(M;\mathbb{R})$ the canonical isomorphism above.
\begin{thm}
\label{thm:3.6} If $f: M \rightarrow N$ is a Lipschitz map then
$f^* \lambda=\lambda f^*$ .
\end{thm}
The last two theorems follow from Lemma \ref{lem:3.4}, \cite[Lemma
3.1]{GKS} and the general theorems of sheaf theory \cite[11.4.6,
11.4.16]{G}.

It turns out that after relaxing the locality requirement, de
Rham Theorem 
above is in general no longer valid for $H^{*}_{\infty}(M)$. In
other words $H^*_{\infty}(M)$ provides a genuinely different (in
fact metric rather than topological) invariant. De Rham theorem
for $H^{*}_{\infty}(M)$ still holds on compact manifolds. This
leads us to ask the following questions.

Do compact manifolds form the maximal class of spaces $M$ for
which the de Rham theorem holds? When is $0 < \dim H^*_{\infty}(M)
< \infty$ for non-compact $M$ ?

\end{document}